\documentclass[12pt]{article} 
\usepackage{amsmath, amssymb} 
\usepackage{amsthm} 
\usepackage{enumerate, url} 
\usepackage{graphicx}

\newtheoremstyle{plainsl}%
	{\topsep}
	{\topsep}
	{\slshape} 
	{}
	{\normalfont\bfseries}
	{.}
	{ }
	{}

\swapnumbers

{\theoremstyle{plainsl}
\newtheorem{theorem}{Theorem}[section]
\newtheorem{lemma}[theorem]{Lemma}
\newtheorem{corollary}[theorem]{Corollary}}
{\theoremstyle{remark}
}

\newcommand\lref[1]{Lemma~\ref{lem:#1}}
\newcommand\tref[1]{Theorem~\ref{thm:#1}}
\newcommand\cref[1]{Corollary~\ref{cor:#1}}

\renewcommand\proof{\noindent\textsl{Proof. }}
\newcommand\sqr[2]{{\vbox{\hrule height.#2pt
    \hbox{\vrule width.#2pt height#1pt \kern#1pt
        \vrule width.#2pt}\hrule height.#2pt}}}
\renewcommand\qed{%
	\ifmmode\eqno\sqr53
	\else\nolinebreak\ \hfill\sqr53\medbreak\fi}

\numberwithin{equation}{section}

\newcommand\be{\beta}
\newcommand\de{\delta}
\newcommand\De{\Delta}

\newcommand\ga{\gamma}

\newcommand\sg{\sigma}

\renewcommand\th{\theta} 

\newcommand\vphi{\varphi}

\newcommand\fld{{\mathbb F}}

\newcommand\ints{{\mathbb Z}}
\newcommand\re{{\mathbb R}}
\newcommand\rats{{\mathbb Q}}

\newcommand\diff{\mathbin{\mkern-1.5mu\setminus\mkern-1.5mu}}

\newcommand\seq[3]{#1_{#2},\ldots,#1_{#3}}

\DeclareMathOperator{\aut}{Aut}

\newcommand\pmat[1]{\begin{pmatrix} #1 \end{pmatrix}}
\newcommand\sm[3]{\sum_{#1=#2}^{#3}}
\DeclareMathOperator{\rk}{rk}
\DeclareMathOperator{\tr}{tr}

\title{When can Perfect State Transfer Occur?} 

\author{
	Chris Godsil\\
	Combinatorics \& Optimization\\
	University of Waterloo\\[2pt]
	\texttt{cgodsil@uwaterloo.ca}}

\begin{document}
\maketitle
	
\begin{abstract}
	Let $X$ be a graph on $n$ vertices with with adjacency matrix $A$ and 
	let $H(t)$  denote the matrix-valued function $\exp(iAt)$.  
	If $u$ and $v$ are distinct
	vertices in $X$, we say \textsl{perfect state transfer} from $u$ to $v$
	occurs if there is a time $\tau$ such that $|H(\tau)_{u,v}|=1$.
	Our chief problem is to characterize the cases where perfect state
	transfer occurs. We show that if perfect state transfer
	does occur in a graph, then the square of its spectral radius is either an integer
	or lies in a quadratic extension of the rationals.
	From this we deduce that for any integer $k$ there only finitely many graphs with
	maximum valency $k$ on which perfect state transfer occurs.
	We also show that if perfect state transfer from $u$ to $v$ occurs,
	then the graphs $X\diff u$ and $X\diff v$ are cospectral and any automorphism
	of $X$ that fixes $u$ must fix $v$ (and conversely).
\end{abstract}

\section{Introduction}

Let $X$ be a graph on $n$ vertices with with adjacency matrix $A$ and 
let $H(t)$  denote the matrix-valued function $\exp(iAt)$.  
If $u$ and $v$ are distinct
vertices in $X$, we say \textsl{perfect state transfer} from $u$ to $v$
occurs if there is a time $\tau$ such that $|H(\tau)_{u,v}|=1$.
The terminology, and much of the motivation, comes from Quantum Physics,
see e.g., \cite{cddekl-pra}. We say that $X$ is periodic relative to a vertex $u$
if there is a time $\tau$ such that $|H(\tau)_{u,u}|=1$, and
say $X$ itself is periodic if there is a time $\tau$ such that
$|H(\tau)_{u,u}|=1$ for all vertices $u$. It is not hard to see
that $X$ will be periodic with period $2\pi$ if all its eigenvalues
are integers. In \cite{cg-period} we showed that $X$ is periodic if and only
if either its eigenvalues are integers, or $X$ is bipartite and there
is a square-free integer $\De$ such that all eigenvalues of $X$
are integer multiples of $\sqrt\De$. As shown by Christandl et al \cite{cddekl-pra}
the $d$-cube is periodic, as are
the Cartesian powers of the path $P_3$. (The eigenvalues
of the $d$-cube are integers, and the eigenvalues of a Cartesian power
of $P_3$ are integer multiples of $\sqrt2$.) 

Our work in this paper has arisen from our attempts to characterize
the pairs of vertices for which perfect state transfer does occur.
Our results taken together imply that finding examples of perfect
state transfer is not going to be easy. In particular we prove that
if perfect state transfer from $u$ to $v$ occurs then
\begin{enumerate}[(1)]
    \item
    Any automorphism of $X$ that fixes $u$ must fix $v$, and conversely.
    \item 
    The subgraphs $X\diff u$ and $X\diff v$ are cospectral.
    \item
    If $\rho$ is the spectral radius of $X$, then $\rho^2$ is an integer.
\end{enumerate}

\section{Periodicity and Perfect State Transfer}

Our first result is elementary, but it shows that perfect state
transfer implies periodicity at at least two vertices. This is
important conceptually, and also plays an important role in the
rest of this paper. If $S$ is a subset of the vertices of the graph $X$
then we use $e_S$ to denote the characteristic vector of $S$. Here the only case
of interest to us will be when $S$ is a single vertex, and then if $u\in V(X)$
then $e_u$ is one of the standard basis vectors.

\begin{lemma}\label{lem:pstper}
	If perfect state transfer from $u$ to $v$ takes place at time $\tau$, then
	it also takes place from $v$ to $u$.  Further $X$ is periodic with
	period $2\tau$ at both $u$ and $v$.
\end{lemma}

\proof
We may assume that
\[
	H(\tau)e_u = \ga e_v
\]
where $|\ga|=1$.  Since $H(t)$ is symmetric for all $t$, we also have
\[
	H(\tau)e_v = \ga e_u
\]
and therefore 
\[
	H(2\tau) e_u =\ga^2 e_u,\quad H(2\tau)e_v =\ga^2 e_v.\qed
\]

A more explicit form of this result is that if 
perfect state transfer from $u$ to $v$ takes place at time $\tau$,
then $H(\tau)$ can be written in block form
\[
    \pmat{H_1&0\\ 0&H_2}
\]
where $H_1$ and $H_2$ are unitary and
\[
    H_1 = \ga \pmat{0&1\\ 1&0}.
\]

We also see that $e_u+e_v$ is an eigenvector for $H(\tau)$
with eigenvalue $\ga$, while $e_u-e_v$ is an eigenvector with
eigenvalue $-\ga$. This has consequences for the spectral
decomposition of $H(\tau)$, which we will discuss after
introducing some associated notation.
If $\seq\th1d$ are the distinct eigenvalues
of $A$, then $E_r$ is the matrix that represents orthogonal projection
onto the eigenspace with eigenvalue $\th_r$. Hence the matrices
$\seq E1d$ are idempotent and symmetric and for any function
$f(x)$ defined on the eigenvalues of $A$, we have
\[
    f(A) = \sum_{r=1}^d f(\th_r)E_r.
\]

It follows that
\[
    H(t) = \sum_r \exp(i\th_r t) E_r;
\]
when $t=\tau$ and we have perfect state transfer from $u$ to $v$
at time $\tau$. We can rewrite this in the form
\[
    H(\tau) = \ga F_+ - \ga F_- +\sum \exp(i\th_r t) E_r.
\]
where the sum is over the integers $r$ such that $\th_r$ is not congruent modulo
$2\pi$ to $\pm\ga$. The matrices $F_+$ and $F_-$ are non-empty sums of the 
principal idempotents and $F_+F_-=0$.

\section{The Ratio Condition}

If $u$ is a vertex in $X$, we define its \textsl{eigenvalue support}
to be the set of eigenvalues $\th_r$ of $X$ such that
\[
    E_r e_u \ne 0.
\]
Note that
\[
    \|E_r e_u\|^2 = e_u^T E_r^TE_r e_u = e_u^T E_r^2 e_u = e_u^T E_re_u
\]
and therefore $\th_r$ lies in the eigenvalue support of $u$ if and only
if $(E_r)_{u,u}\ne0$.
The ratio condition is a necessary condition for a graph to be periodic
at a vertex. The version we offer here is stated as Theorem~2.2 in \cite{cg-period};
it is an extension of result from Saxena, Severini and Shparlinski \cite{sasesh},
which in turn extends an idea used in Christandl et al \cite{cddekl-pra}.

\begin{theorem}\label{thm:rat}
    Let $X$ be a graph and let $u$ be a vertex in $X$ at which $X$ is periodic.
    If $\th_k$, $\th_\ell$, $\th_r$, $\th_s$ are eigenvalues in the support
    of $e_u$ and $\th_r\ne\th_s$, then
    \[
        \frac{\th_k-\th_\ell}{\th_r-\th_s} \in \rats.\qed
    \]
\end{theorem}

\section{Orbits and Equitable Partitions}

If $u\in V(X)$, then $\aut(X)_v$ denotes the group of automorphisms of $X$ 
that fix $v$.  Our next result says that if perfect state transfer from $u$ 
to $v$ occurs, then any automorphism of $X$ that fixes $u$ must fix $v$ 
(and vice versa).

\begin{lemma}\label{lem:pst-stab}
	If $X$ admits perfect state transfer from $u$ to $v$, then 
	$\aut(X)_u=\aut(X)_v$.
\end{lemma}

\proof
We identify the automorphisms of $X$ with the permutation matrices that 
commute with $A$. Since $H(t)$ is a polynomial in $A$, any permutation matrix 
from $\aut(X)$ commutes with $H(\tau)$.
If the automorphism associated with $P$ fixes $u$, the $Pe_u=e_u$.
If perfect state transfer takes place at time $\tau$ and $H=H(\tau)$, then there
is a complex scalar $\ga$ such that $|\ga|=1$ and $He_u=\ga e_v$.  So
\[
	\ga Pe_v = PHe_u = HPe_u =He_u =\ga e_v
\]
and thus $v$ is fixed.  As our argument is symmetric in $u$ and $v$, 
the lemma follows.\qed

By way of example, suppose $X$ is a Cayley graph for an abelian
group $G$. The map that sends $g$ in $G$ to $g^{-1}$ is an automorphism
of $G$ and also an automorphism of $X$. The fixed points of this
automorphism are the elements of $G$ with order one or two.
So if perfect state transfer for 1 to a second vertex $a$ occurs,
then $a$ has order two. Consequently perfect state transfer cannot
occur on a Cayley graph for an abelian group of odd order. (Another
proof of this is presented as Corollary~4.2 in \cite{cg-period}.)

A partition $\pi$ of the vertices of a graph $X$ is \textsl{equitable}
if, for each ordered pair of cells $C_i$ and $C_j$ from $\pi$, all vertices 
in $C_i$ have the same number of neighbors in $C_j$. (So the subgraph of
$X$ induced by a cell is a regular graph, and the subgraph formed by the
vertices of two cells and the edges which join them is bipartite and
semiregular.) We note that the orbits of any group of automorphisms of $X$ 
form an equitable partition. The \textsl{discrete partition}, with each cell 
a singleton,
is always equitable; the \textsl{trivial partition} with exactly one cell
is equitable if an only if $X$ is regular. For the basics concerning equitable 
partitions see for example \cite[Sction~9.3]{cggrbk}. The 
join of two equitable partitions is equitable, and consequently given any 
partition of $X$, there is a unique coarsest refinement of it---the join of 
all equitable partitions which refine it. 

If $u\in V(G)$, we use $\De_u$ to denote the coarsest equitable
refinement of the partition
\[
    \{ \{u\}, V(X)\diff\{u\} \}
\]
The orbits of the group of automorphisms of $X$ that fix $u$ will be
a refinement of $\De_u$. 

If $\pi$ is a partition of $V(X)$, its \textsl{characteristic matrix} 
is the $01$-matrix
whose columns are the characteristic vectors of the cells of $\pi$, viewed
as subsets of $V(X)$. If we scale the columns of the characteristic matrix
so they are unit vectors, we obtain the \textsl{normalized characteristic matrix}
of $\pi$. If $P$ and $Q$ are respectively the characteristic and normalized
characteristic matrix of $\pi$, then $P$ and $Q$ have the same column space 
and $Q^TQ=I$. The matrix $QQ^T$ is block diagonal, and its diagonal blocks are
all of the form
\[
    \frac1r J_r
\]
where $J_r$ is the all-ones matrix of order $r\times r$, and the size 
of the $i$-th block is the size of the $i$-th cell of $\pi$.
The vertex $u$ forms a singleton cell of $\pi$ if and only if $Qe_u=e_u$.

We use $|\pi|$ to denote the number of cells of $\pi$.

\begin{lemma}
    Suppose $\pi$ is a partition of the vertices of the graph $X$,
    and that $Q$ is its normalized characteristic matrix. Then
    the following are equivalent:
    \begin{enumerate}[(a)]
        \item 
        $\pi$ is equitable.
        \item
        The column space of $Q$ is $A$-invariant.
        \item
        There is a matrix $B$ of order $|\pi|\times|\pi|$ such that $AQ=QB$.
        \item
        $A$ and $QQ^T$ commute.\qed
    \end{enumerate}
\end{lemma}

If there is perfect state transfer form $u$ to $v$, then it follows from 
\lref{pst-stab} that the orbit partitions of $\aut(X)_u$ and $\aut(X)_v$
are equal. Since equitable partitions are an analog of orbit partitions,
it is not entirely unreasonable to view the following result as an extension
of this fact.

\begin{theorem}
    Let $u$ and $v$ be vertices in $X$. If there is perfect state transfer
    from $u$ to $v$, then $\De_u=\De_v$.
\end{theorem}

\proof
Let $Q$ be the normalized characteristic matrix of the partition
$\De_u$ and assume $H(t)=\exp(iAt)$. Then $H(t)$ is a polynomial
and so $QH(t)=H(t)Q$. If we have perfect state transfer from $u$
to $v$ at time $\tau$, then
\[
    H(\tau)e_u = \ga e_v
\]
and accordingly
\[
    \ga Qe_v = QH(\tau)e_u = H(\tau) Qe_u = H(\tau)e_u = \ga e_v.
\]
So $Qe_v=e_v$ and this implies that $\{v\}$ is a cell of $\De_u$.\qed

By way of example, consider a distance-regular graph $X$ with diameter $d$.
If $u\in V(X)$, then $\De_u$ is the distance partition
with respect to $u$ with exactly $d+1$ cells, where the $i$-th
cell is the set of vertices at distance $i$ from $u$.
We conclude that if perfect state transfer occurs on $X$, then for each vertex
$u$ there is exactly one vertex vertex at distance $d$ from it.
In particular perfect state transfer does not occur on strongly regular graphs.
(This also follows from some observations in \cite[Section~4]{cg-period}.)

\section{Cospectral Vertices}

Vertices $u$ and $v$ in the graph $X$ are \textsl{cospectral} if
\[
	\phi(X\diff u,t) = \phi(X\diff v,t).
\]
Of course two vertices that lie in the same orbit of $\aut(X)$ are cospectral,
but there are many examples of cospectral pairs of vertices where this does 
not hold. A graph is \textsl{walk regular} if any two of its vertices 
are cospectral. Any strongly regular graph is walk regular.

\begin{lemma}\label{lem:pst-cosp}
	If $X$ admits perfect state transfer from $u$ to $v$, then 
	$E_ru=\pm E_rv$ for all $r$, and $u$ and $v$ are cospectral.
\end{lemma}

\proof
By Cramer's rule,
\[
    ((tI-A)^{-1})_{u,u} = \frac{\phi(X\diff u,t)}{\phi(X,t)}
\]
and so using the spectral decomposition we find that
\[
	\frac{\phi(X\diff u,t)}{\phi(X,t)} =\sum_r (t-\th_r)^{-1} (E_r)_{u,u},
\]
where
\[
	(E_r)_{u,u} = \| E_r e_u \|^2.
\]
Assume we have perfect state transfer from $u$ to $v$ at time $\tau$ 
and set $H=H(\tau)$. Then there is a complex number $\ga$ such that 
$|\ga|=1$ and $He_u =\ga e_v$.

As $HE_r = \exp(i\th_r\tau)E_r$ we have
\[
    \exp(i\th_r\tau)E_r e_u = \ga E_re_v
\]
and since $E_re_u$ and $E_re_v$ are real our first claim follows.
This implies that
\[
    e_u^TE_re_u = \| E_re_u\|^2 = \| E_re_v\|^2 = e_v^TE_re_v
\]
and therefore $u$ and $v$ are cospectral.\qed

\section{Integrality, Nearly}

From \cite{cg-period} we know that if a graph is periodic, then the squares of 
its eigenvalues are integers and, if the graph is not bipartite, 
the eigenvalues themselves are integers. Here we derive a variant of this
fact, which will imply that that if perfect state transfer occurs
on $X$, then the square of the spectral radius of $X$ is an integer.

At the end of the previous section we saw
that if we had perfect state transfer from $u$ to $v$ at time $\tau$,
then
\[
    H(\tau)E_r e_u = \ga E_r e_v
\]
from which it follows that the eigenvalue supports of $u$ and $v$ are equal.

\begin{theorem}
    \label{thm:intevs}
    Suppose $X$ is a connected graph with at least three vertices
    where perfect state transfer
    from $u$ to $v$ occurs at time $\tau$. Let $\de$ be the dimension of the
    $A$-invariant subspace generated by $e_u$. Then $\de\ge3$ and either
    all eigenvalues in the eigenvalue support of $u$ are integers, or
    they are all of the form $\frac12(a+b_i\sqrt\De$) 
    where $\De$ is a square-free integer and $a$ and $b_i$ are integers.
\end{theorem}

\proof
Let $S$ denote the eigenvalue support of $e_u$. From the spectral decomposition,
the $A$-invariant subspace generated by $e_u$ is spanned by the non-zero
vectors $E_re_u$ such that $\th_r\in S$. Since $E_rE_s=0$ when $r\ne s$
these vectors are pairwise orthogonal and so form a basis for our $A$-invariant
subspace. We conclude that $\de=|S|$.

Since we have transfer from $u$ to $v$ at time $\tau$, there is a complex 
number $\ga$ such that
\[
    H(\tau)e_u = \ga e_v
\]
from which we see that that the $A$-invariant subspace generated by $e_u$ contains
$e_v$. Hence $\de\ge2$, but if equality holds then the span of
$e_u$ and $e_v$ is $A$-invariant and, since $X$ is connected it
follows that $|V(X)|=2$.

If $\th_r$ and $\th_s$ are algebraic conjugates, then
$E_r$ and $E_s$ are algebraic conjugates and so $E_re_u=0$ if and only if $E_se_u=0$.
Hence $S$ contains all algebraic conjugates of each of its elements.
If two eigenvalues in $S$ are integers, say $\th_0$ and $\th_1$ then
since the ratio conditions asserts that
\[
    \frac{\th_r-\th_0}{\th_1-\th_0} \in \rats,
\]
we conclude that all elements of $S$ are integers. 

If all eigenvalues in $S$ are integers, we're done, so suppose
that $\th_0$ and $\th_1$ are two distinct irrational eigenvalues in $S$. 
We show that $(\th_1-\th_0)^2$ is an integer. By the ratio condition, 
if $\th_r,\th_s\in S$ there is a rational number $a_{r,s}$ such that
\[
    \th_r-\th_s = a_{r,s}(\th_1-\th_0)
\]
and therefore
\[
    \prod_{r\ne s}(\th_r-\th_s) 
        =(\th_1-\th_0)^{\de^2-\de}\prod_{i\ne j}a_{r,s}.
\]
The product on the left is an integer and the product of the 
$a_{r,s}$'s is rational, and hence
\[
    (\th_1-\th_0)^{\de^2-\de} \in \rats.
\]
Since $\th_1-\th_0$ is an algebraic integer, this implies that
\[
    (\th_1-\th_0)^{\de^2-\de} \in \ints.
\]
Suppose $m$ is the least positive integer such that $(\th_1-\th_0)^m$ 
is an integer. Then there are $m$ distinct conjugates of $\th_1-\th_0$ of the form
\[
    \be e^{2\pi ik/m}\qquad (k=0,\ldots,m-1)
\]
where $\be$ is the positive real $m$-th root of an integer,
and since the eigenvalues of $X$ are real we conclude that $m\le2$.
Since $\th_1-\th_0$ is not rational, it must therefore be the square root 
of an integer. Assume $\De=(\th_1-\th_0)^2$.

Let $\fld$ be the extension of $\rats$ generated by $S$ and let $G$ be the Galois
group of this extension. We aim to show that $|\fld:\rats|=2$. Let $\sg$ be a 
non-identity element of $G$ and let $\th$ be an irrational element in $S$.
Then
\[
     \frac{\th^{\sg^2}-\th^\sg}{\th^\sg-\th} = a \in \rats
\]
and consequently
\[
    \frac{\th^{\sg^{k+1}}-\th^{\sg^k}}{\th^{\sg^k}-\th^{\sg^{k-1}}} 
        = a^{\sg^{k-1}} = a.
\]
The product of all the possible fractions on the left here is 1, so if
$m$ is the order of $\sg$, then $a^m=1$. Since $a$ is real, this implies
that $a=\pm1$. 

If $a=1$, then
\[
    \th^{\sg^{k+1}}-\th^{\sg^k} = \th^{\sg^k}-\th^{\sg^{k-1}},
\]
from which it follows that for all $k$
\[
    \th^{\sg^{k+1}}-\th = k(\th^\sg-\th).
\]
Setting $k$ equal to the order of $\sg$ yields a contradiction, and thus
$a=-1$. Now we have
\[
    \th^{\sg^2} - \th^\sg = -\th^\sg + \th
\]
and accordingly $\sg^2$ fixes $\th$. Therefore all elements of the Galois
group $G$ have order dividing two. 

Now let $\vphi$ be an element of $S$. Since
\[
    \frac{\vphi-\th}{\th^\sg-\th} \in \rats
\]
we have
\[
    \vphi = a\th^\sg + (1-a)\th
\]
for some rational $a$. If $\rho\in G$, then there is a rational $b$ such that
\begin{align*}
    \th^\rho &= b \th^\sg + (1-b)\th\\
    \th^{\sg\rho} &= (1-b)\th^\sg +b\th.
\end{align*}
As $\rho^2=1$ and
\[
    \pmat{b&1-b\\ 1-b&b}^2 = \pmat{b^2+(1-b)^2& 2b(1-b)\\ 2b(1-b)& b^2+(1-b)^2}
\]
we have $b-b^2=0$. If $b=0$ then $\rho=1$, if $b=1$ then $\rho=\sg$.
We conclude that $|G|=2$ and $|\fld:\rats|=2$.

As $|G|=2$, the field $\rats(\th_0)$ is a quadratic extension of $\rats$
and consequently $\rats(\th_0)=\rats(\sqrt\De)$. Because $X$ is connected
its spectral radius is simple and the corresponding eigenvector is positive.
Therefore the spectral radius must lie in the eigenvalue support $S$ of $e_u$.
Assume $\th_0$ is the spectral radius and suppose it is an integer.
Then if $\th_1$ is an irrational eigenvalue in $S$, we may assume that
\[
    \th_1 = \th_0 -\sqrt{\De}
\]
(where $\De\in\ints$) but then the conjugate $\th_0+\sqrt{\De}$ is an eigenvalue
of $X$ greater than $\th_0$. Therefore if $\th_0$ is an integer, all elements
of $S$ are integers.

Assume $\th_0$ is irrational, and let $\seq\th0m$ be the irrational eigenvalues
in $S$. Then $\th_0-\th_1=\sqrt\De$ for some integer $\De$ and for $i=0,\ldots,m$
there are rationals $a_i$ and $b_i$ such that
\[
    \th_i = a_i+b_i\sqrt\De.
\] 
(Since $\th_i$ is an algebraic integer, both $2a_i$ and $2b_i$ must be integers.)
The conjugate $\bar\th_i$ of $\th_i$ is $a_i-b_i\sqrt\De$ and, 
by the ratio condition,
\[
    \frac{\th_i-\th_0}{\bar\th_0-\th_0} =
        \frac{a_i-a_0+(b_i-b_0)\sqrt\De}{2b_0\sqrt\De} \in\rats,
\]
whence $a_i=a_0$ for $i=0,\ldots,m$.

If there is an integer eigenvalue $\psi$ in $S$, then
\[
    \frac{\psi-a_0-b_i\sqrt\De}{\psi-a_0+b_i\sqrt\De} \in \rats.
\]
Since
\[
    \frac{a+b\sqrt\De}{c+d\sqrt\De}\in\rats \iff \det\pmat{a&b\\c&d} = 0,
\]
it follows that $\psi=a_0$.\qed

Suppose $\de=3$. Then the subspace $U$ generated by $e_u$ contains $e_v$, $Ae_u$
and $Ae_v$. If $Ae_u=e_v$, then $\de=2$, and so $U$ is spanned by $e_u$,
$e_v$ and either one of $Ae_u$ and $Ae_v$. We find that $X$ is the join
of $u$ and $v$ with a regular graph.

\begin{corollary}
    There are only finitely many connected graphs with maximum valency at most $k$
    where perfect state transfer occurs.
\end{corollary}

\proof
Suppose $X$ is a connected graph where perfect state transfer
from $u$ to $v$ occurs at time $\tau$ and let $S$ be the eigenvalue support
of $u$. If the eigenvalues in $S$ are integers then $|S|\le2k+1$, and if they
are not integers then $|S|<(2k+1)\sqrt{2}$. So the dimension of the $A$-invariant
subspace of $\re^{V(X)}$ generated by $e_u$ is at most $\lceil(2k+1)\sqrt{2}\rceil$, 
and this is also a bound on the maximum distance from $u$ of a vertex in $X$.
If $s\ge1$, the number of vertices at distance $s$ from $u$ is at most
$k(k-1)^{s-1}$, and the result follows.\qed

\section{Controllable Vertices}

Let $A$ be the adjacency matrix of the graph $X$ on $n$ vertices and suppose 
$S$ is a subset of $V(X)$ with characteristic vector $e_S$. The \textsl{walk matrix}
of $S$ is the matrix with columns
\[
    e_S, Ae_S,\ldots, A^{n-1}e_S.
\]
We say that the pair $(X,S)$ is \textsl{controllable} if this walk matrix is
invertible. If $u\in V(X)$ then the walk matrix of $u$ is just the walk matrix
relative to the subset $\{u\}$ of $V(X)$. From our discussion at the start of the
proof of \tref{intevs},
we have that the rank of the walk matrix relative to $u$ is equal to the
size of the eigenvalue support of $u$. Thus if $(X,u)$ is controllable,
the eigenvalues of $A$ must be distinct.

We say that vertices $u$ and $v$ of $X$ are \textsl{cospectral} if the
vertex-deleted subgraphs
\[
    X\diff u,\quad X\diff v
\]
are cospectral. We use $\phi(X,x)$ to denote the characteristic polynomial
of $A(X)$.

\begin{lemma}\label{lem:wutwu}
    Let $u$ and $v$ be vertices of $X$ with respective walk matrices $W_u$ 
    and $W_v$. Then $u$ and $v$ are cospectral if and only if 
    $W_u^TW_u=W_v^TW_v$.
\end{lemma}

\proof
By Cramer's rule
\[
    ((xI-A)^{-1})_{u,u} = \frac{\phi(X\diff u,x)}{\phi(X,x)}
\]
Now
\[
    ((xI-A)^{-1})_{u,u} = x^{-1}\sum_{k\ge0} (A^k)_{u,u} x^{-k}
\]
and therefore vertices $u$ and $v$ are cospectral if and only if
\[
    (A^k)_{u,u} =(A^k)_{v,v}
\]
for all non-negative integers $k$. If $n=|V(X)|$, then equality holds
for all non-negative integers $k$ if and only if it holds for
$k=0,1,\ldots,n-1$. To complete the proof observe that if $W_u$ is the
walk matrix for the vertex $u$, then
\[
    (W_u^TW_u)_{r,s} = e_u^TA^{r+s-2}e_u.\qed
\]

One consequence of this is that if $u$ and $v$ are cospectral and
$(X,u)$ is controllable, then so is $(X,v)$. We also have:

\begin{corollary}\label{cor:}
    If $u$ is a vertex in $X$ with walk matrix $W_u$, then $\rk(W_u)$
    is equal to the number of poles of the rational function
    $\phi(X\diff u,x)/\phi(X,x)$. Hence $(X,u)$ is controllable if and only
    if $\phi(X\diff u,x)$ and $\phi(X,x)$ are coprime.
\end{corollary}

\proof
By the spectral decomposition,
\[
    ((xI-A)^{-1})_{u,u} = \sum_r \frac{1}{x-\th_r}(E_r)_{u,u}.\qedß
\]

We will prove that a controllable vertex cannot be involved in perfect state
transfer. For this we need the following:

\begin{lemma}\label{lem:evsep}
    Let $X$ be a graph with $n=|V(X)|$.
    If $\sg$ is the minimum distance between two eigenvalues of $X$, then
    \[
        \sg^2 < \frac{12}{n+1}.
    \]
\end{lemma}

\proof
Assume that the eigenvalues of $X$ in non-increasing order are $\seq\th1n$.
If we have
\[
    M := A\otimes I - I\otimes A
\]
then the eigenvalues of $M$ are the numbers
\[
    \th_i-\th_j,\quad 1\le i,j \le n.
\]
Now
\[
    M^2 =A^2\otimes I +I\otimes A^2 -2A\otimes A
\]
and consequently, if $e:=|E(X)|$ then
\[
    \sm{i,j}1n(\th_i-\th_j)^2 = \tr(M^2) =2n\tr(A^2) = 4ne.
\]

Since
\[
    \th_i-\th_j \ge (i-j)\sg
\]
we have
\[
    \sm{i,j}1n (\th_i-\th_j)^2 \ge \sg^2 \sm{i,j}1n (i-j)^2
\]
As 
\[
    \sm i1n i^2 = \frac{n(n+1)(2n+1)}{6}
\] 
we find that
\[
    \sm{i,j}1n (i-j)^2 =2n\sm i1n i^2 -2\left(\sm i1n i\right)^2
        =\frac{n^2(n^2-1)}{6},
\]
and since $e\le n(n-1)$ this yields
\[
    \sg^2\frac{n^2(n^2-1)}{6} \le 4ne \le 2n^2(n-1).
\]
This gives our stated bound but with $\le$ in place of $<$. To achieve
strictness we note that if equality were to hold then $e=n(n-1)$
and $X=K_n$. Since $\sg(K_n)=0$, we are done.\qed

\begin{theorem}\label{thm:nocontrol}
    Let $X$ be a connected graph on at least four vertices.
    If we have perfect state transfer between distinct vertices $u$ and 
    $v$ in $X$, then neither $u$ nor $v$ is controllable.
\end{theorem}

\proof
Let $n=|V(X)|$. We assume by way of contradiction that $(X,u)$ is controllable.
Then the eigenvalue support of $u$ contains all eigenvalues of $X$
and these eigenvalues are distinct. By \tref{intevs}, this means that
there are integers $a$ and $\De$ and distinct integers $\seq b1n$
and the eigenvalues of $A$ are the numbers 
\[
    \frac12(a+b_r\sqrt{\De}).
\]
As $\tr(A)=0$ we have
\[
    0 = \frac12 \sum_{r=1}^n (a+b_i\sqrt{\De}) 
        = \frac{na}{2} +\frac{\sqrt{\De}}{2}\sum_{r=1}^n b_r
\]
If $\sqrt{\De}$ is not an integer, this implies that $\sum_r b_r=0$
and hence that $a=0$. Accordingly the eigenvalues of $A$ are distinct
integer multiples of the irrational number $\sqrt{\De}$ that sum to zero. If 
$\sqrt{\De}$ is an integer then the eigenvalues of $A$ are distinct integers
that sum to zero. In either case $\sg(X)\ge1$.

By \lref{evsep} we can assume that $n=|V(X)\le 10$. This leaves us with six cases. 
First suppose $\De=1$. Assume $n=10$. Then the sum of the squares of the
eigenvalues of $X$ is bounded below by the sum of the squares of the integers 
from $-4$ to $5$, which is 85, and hence the average valency of a vertex is 
at least $8.5$. 
This implies that $\th_1=9$, not $5$, consequently the sum of the squares 
of the eigenvalues is at least
\[
    85 -25 +81 = 141
\]
and now the average valency is $14.1$, which is impossible.
The cases where $n$ is 7, 8 or 9 all yield contradictions in the same way.
If $\De>1$, it is even easier to derive contradictions.

Next, brute force computation (using Sage \cite{sage}) shows that the path $P_4$
is the only graph on 4, 5 or 6 vertices where the minimum separation between 
consecutive eigenvalues is at least 1. The positive eigenvalues of $P_4$ are
\[
    (\sqrt{5}\pm1)/2
\]
and their ratio is not rational.\qed

\section{Nearly Perfect State Transfer?}

If $u$ and $v$ are controllable and cospectral we cannot get perfect state transfer
between them. Our next results shows that if these conditions hold, there is
a symmetric orthogonal matrix $Q$ which commutes with $A$ and maps $e_u$ to $e_v$.

\begin{lemma}\label{lem:}
    Let $u$ and $v$ be vertices in $X$ with respective walk matrices $W_u$ and $W_v$.
    If $u$ and $v$ are controllable and $Q:=W_vW_u^{-1}$ then $Q$ is a polynomial in $A$.
    Further $Q$ is orthogonal if and only if $u$ and $v$ are cospectral.
\end{lemma}

\proof
Let $C_\phi$ denote the companion matrix of the characteristic polynomial of $A$. Then
\[
    AW_u =W_u C_\phi
\]
for any vertex $u$ in $X$. Hence if $u$ and $v$ are controllable,
\[
    W_u^{-1}AW_u = W_v^{-1}AW_v
\]
and from this we get that
\[
    AW_vW_u^{-1} = W_vW_u^{-1}A.
\]
Since $X$ has a controllable vertex its eigenvalues are all simple, and so any matrix
that commutes with $A$ is a polynomial in $A$. This proves the first claim. 

From \lref{wutwu}, the vertices $u$ and $v$ are cospectral
if and only if $W_u^TW_u=W_v^TW_v$, which is equivalent to
\[
    Q^{-T} = W_v^{-T}W_u^T = W_vW_u^{-1} = Q.\qed
\]

This places us in an interesting position. If $u$ and $v$ are cospectral and
controllable, there is an orthogonal matrix $Q$ that commutes with $A$ such that
$Qe_u=e_v$ but, by the result of the previous section, $Q$ cannot be equal to
a scalar multiple of $H(t)$ for any $t$.

\section*{Acknowledgement}

I thanks Dave Morris for carefully reading the first draft of this paper and
for providing a number of useful comments.


\begin{thebibliography}{1}

\bibitem{cddekl-pra}
{\sc M.~Christandl, N.~Datta, T.~Dorlas, A.~Ekert, A.~Kay, and A.~Landahl},
  {\em Perfect transfer of arbitrary states in quantum spin networks}, Phys Rev
  A, 71 (2005), p.~032312.

\bibitem{cg-period}
{\sc C.~Godsil}, {\em Periodic graphs}, arXiv, math.CO (2008), 0806.2074v2.
\newblock 19 pages.

\bibitem{cggrbk}
{\sc C.~Godsil and G.~Royle}, {\em Algebraic {G}raph {T}heory}, vol.~207 of
  Graduate Texts in Mathematics, Springer-Verlag, New York, 2001.

\bibitem{sasesh}
{\sc N.~Saxena, S.~Severini, and I.~Shparlinski}, {\em Parameters of integral
  circulant graphs and periodic quantum dynamics}, arXiv, quant-ph (2007),
  quant-ph/0703236v1.
\newblock 12 pages.

\bibitem{sage}
{\sc W.~Stein et~al.}, {\em {S}age {M}athematics {S}oftware ({V}ersion 4.5.2)},
  The Sage Development Team, 2010.
\newblock {\tt http://www.sagemath.org}.

\end{thebibliography}

\end{document}